\documentclass[12]{article}

\usepackage{latexsym}
\usepackage{amsfonts}
\usepackage{graphicx}
\usepackage{amsmath}
\usepackage{amssymb}
\usepackage[latin1]{inputenc}

\title{\textbf{On a relation of bicomplex pseudoanalytic function theory to the complexified stationary Schr{\"o}dinger equation}}
\author{D. Rochon\thanks{Research supported by CRSNG(Canada) and FQRNT(Qu{\'e}bec)}}
\date{D\'epartement de math\'ematiques et
d'informatique \\ Universit\'e du Qu\'ebec \`a Trois-Rivi\`eres \\
C.P. 500 Trois-Rivi\`eres, Qu\'ebec \\ Canada, G9A 5H7
\\E-mail: Dominic.Rochon@UQTR.CA}

\makeatletter \@addtoreset{equation}{section}

\makeatletter

\def\be   {\begin{equation}}   \def\ee   {\end{equation}}
\def\ba   {\begin{array}}      \def\ea   {\end{array}}
\def\bea  {\begin{eqnarray}}   \def\eea  {\end{eqnarray}}
\def\bean {\begin{eqnarray*}}  \def\eean {\end{eqnarray*}}

\newtheorem{theorem} {Theorem}
\newtheorem{lemma}{Lemma}
\newtheorem{definition} {Definition}

\newtheorem{remark}{Remark}

\newtheorem{proposition} {Proposition}

\newcommand{\pre}{\mathrm{Re}}
\newcommand{\pim}{\mathrm{Im}}

\newcommand{\bi} {\ensuremath{{\bf i}}}
\newcommand{\bo} {\ensuremath{{\bf i_1}}}
\newcommand{\bos}{\ensuremath{{\bf i_1^{\text 2}}}}
\newcommand{\bt} {\ensuremath{{\bf i_2}}}
\newcommand{\bts}{\ensuremath{{\bf i_2^{\text 2}}}}

\newcommand {\bj}{\ensuremath{{\bf j}}}
\newcommand {\bje}{\ensuremath{{\bf j^{*}}}}
\newcommand {\bjs}{\ensuremath{{\bf j^{\text 2}}}}
\newcommand{\eo} {\ensuremath{{\bf e_1}}}
\newcommand{\et} {\ensuremath{{\bf e_2}}}

\newcommand{\mC}{\ensuremath{\mathbb{C}}}

\newcommand{\mR}{\ensuremath{\mathbb{R}}}

\begin{document}
\maketitle

\begin{abstract}
Using three different representations of the bicomplex numbers $\mathbb{T}\cong {\rm Cl}_{\Bbb{C}}(1,0) \cong {\rm Cl}_{\Bbb{C}}(0,1)$,
which is a commutative ring with zero divisors defined by
$\mathbb{T}=\{w_0+w_1 {\bf i_1}+w_2{\bf i_2}+w_3 {\bf j}\ |\ w_0,w_1,w_2,w_3 \in
\mathbb{R}\}$ where ${\bf i_1^{\text 2}}=-1,\ {\bf i_2^{\text 2}}=-1,\ {\bf j^{\text 2}}=1 \mbox{ and }\ {\bf i_1}{\bf i_2}={\bf j}={\bf i_2}{\bf i_1}$,
we construct three classes of bicomplex pseudoanalytic functions. In particular, we obtain
some specific systems of Vekua equations of two complex variables and we established some connections between one of these systems
and the classical Vekua equations. We consider also the complexification of the real stationary two-dimensional Schr{\"o}dinger equation. With the
aid of any of its particular solutions, we construct a specific bicomplex Vekua equation possessing the following special property.
The scalar parts of its solutions are solutions of the original complexified Schr{\"o}dinger equation and the vectorial parts are solutions
of another complexified Schr{\"o}dinger equation.
\end{abstract}

\vspace{1cm}

\normalsize
\newpage

\section{Introduction}

The pseudoanalytic function theory was independently developed by two prominent mathematicians, I.N. Vekua (see \cite{17}) and L. Bers (see \cite{1,3,4}).
Historically, the theory became one of the important impulses for developing the general theory of elliptic systems.
More recently, it has been established by V.V. Kravchenko (see \cite{8,9}) that with the aid of any particular solutions of the real stationary two-dimensional
Schr{\"o}dinger equation we can construct a Vekua equation possessing the following special property. The real parts of its solutions are solutions of the original Schr{\"o}dinger equation
and the imaginary parts are solution of an associated Schr{\"o}dinger equation with a potential having the form of a potential obtained after
a Darboux transform. Using Bers's theory of Taylor series for pseudoanalytic function, the author obtain a locally complete system of solutions
of the original Schr{\"o}dinger equation which can be constructed explicitly for an ample class of Schr{\"o}dinger equation. Subsequently, V.V. Kravchenko (see \cite{19})
gives a generalization of the factorization technics developed in \cite{8} for the more general two-dimensional elliptic operator $(\mbox{div}\mbox{ }p\mbox{ } \mbox{grad} +q)u=0$, and consider
the case where $p$ and $q$ are complex functions. In particular, the author had to consider a bicomplex Vekua equation of a special form to be able
to obtain the following more general property. The scalar parts of the bicomplex Vekua equation considered are solutions of the original Schr{\"o}dinger equation with a complex-valued potential.
However, the case using the complex functions is more complicated and the author let as an open question the proof of expansion and convergence theorems for
the bicomplex Vekua equation considered.

In this article, using three different representations of the bicomplex numbers (see, e.g. \cite{6,11,12,13,14,15,19}),
which is a commutative ring with zero divisors defined by
$\mathbb{T}=\{w_0+w_1 {\bf i_1}+w_2{\bf i_2}+w_3 {\bf j}\ |\ w_0,w_1,w_2,w_3 \in
\mathbb{R}\}$ where ${\bf i_1^{\text 2}}=-1,\ {\bf i_2^{\text 2}}=-1,\ {\bf j^{\text 2}}=1 \mbox{ and }\ {\bf i_1}{\bf i_2}={\bf j}={\bf i_2}{\bf i_1}$,
we construct three classes of bicomplex pseudoanalytic functions. For every class we obtain a
special type of bicomplex Vekua equation of two complex variables, of which the one considered by  A. Castaneda and V.V. Kravchenko (see \cite{5,20})
when the domain is restricted to the complex (in $\bt$) plane. Moreover, we established some connections between one of these systems
of bicomplex Vekua equation and the classical Vekua equations. We consider also the complexification of the real stationary
two-dimensional Schr{\"o}dinger equation : $$(\triangle_{\mathbb{C}}-\nu(z_1,z_2))f=0$$
where $\omega=z_1+z_2\bt\in\mathbb{T}$ with $z_1,z_2\in\mathbb{C}(\bo)$ and $\triangle_{\mathbb{C}}$ is the complexified Laplacian operator i.e. $\triangle_{\mathbb{C}}={\partial^{2}_{z_1}}+{\partial^{2}_{z_2}}$.
With the aid of any of its particular solutions $f_0$ and the bicomplex operators : $$\partial _{{\omega}^{\dag_{2}}} = \frac{1}{2}\left( {\partial _{z_1}  + \bt\partial _{z_2} }\right)
\mbox{ and } \partial _{\omega } = \frac{1}{2}\left( {\partial _{z_1}  - \bt\partial _{z_2} } \right),$$ we obtain the following factorization of the
complexified Schr{\"o}dinger equation
$$(\triangle_{\mathbb{C}}-\nu)\varphi=4\left(\partial_{{\omega}^{\dagger_2}}+\frac{\partial_{\omega} f_0}{f_0}C\right)
\left(\partial_{\omega}-\frac{\partial_{\omega} f_0}{f_0}C\right)\varphi$$
where $C$ denote the $\dagger_2$-bicomplex conjugation operator,
and we consider a specific bicomplex Vekua equation
$$\left( \partial_{{\omega}^{\dagger_2}} - \frac{\partial_{{\omega}^{\dagger_2}}f_0}{f_0}C \right) W=0$$
possessing the following special property.
The scalar parts of its solutions are solutions of the original complexified Schr{\"o}dinger equation
and the vectorial parts are solutions
of another complexified Schr{\"o}dinger equation with the following potential $-\nu(z_1,z_2)+\frac{2|\bigtriangledown_{\mathbb{C}}f_0|^{2}_{\bo}}{f^2_2}$
where $|\omega|^{2}_{\bo}=\omega{\omega}^{\dagger_2}$ $\forall \omega\in\mathbb{T}$ and
$\bigtriangledown_{\mathbb{C}}$ is the complexified gradient operator i.e. $\bigtriangledown_{\mathbb{C}}={\partial_{z_1}}+\bt {\partial_{z_2}}$.

Finally, from the fact that the complexified Schr{\"o}dinger equation contains
the stationary two-dimensional Schr{\"o}dinger equation
$$(\triangle-\nu(x,p))f=0$$
and the Klein-Gordon equation
$$(\square-\nu(x,q))f=0,$$
we show that our factorization of the complexified Schr{\"o}dinger equation is a generalization of
the factorization obtained in \cite{8} for the stationary two-dimensional Schr{\"o}dinger equation and
for the factorization obtained in \cite{24} for the Klein-Gordon equation.


\section{Preliminaries}

\subsection{Bicomplex Numbers}
Bicomplex numbers are defined as

\be \mathbb{T}:=\{z_1+z_2\bt\ |\ z_1, z_2 \in \mathbb{C}(\bo) \}
\label{enstetra}
\ee where the imaginary units $\bo, \bt$ and $\bj$
are governed by the rules: $\bos=\bts=-1$, $\bjs=1$ and \be
\ba{rclrcl}
   \bo\bt &=& \bt\bo &=& \bj,  \\
   \bo\bj &=& \bj\bo &=& -\bt, \\
   \bt\bj &=& \bj\bt &=& -\bo. \\
\ea \ee Note that we define $\mC(\bi_k):=\{x+y\bi_k\ |\ \bi_k^2= -1$
and $x,y\in \mR \}$ for $k=1,2$. Hence, it is easy to see that the
multiplication of two bicomplex numbers is commutative. In fact, the bicomplex numbers
$$\mathbb{T}\cong {\rm Cl}_{\Bbb{C}}(1,0) \cong {\rm Cl}_{\Bbb{C}}(0,1)$$
are \textit{unique} among the complex Clifford algebras
in that they are commutative but not division algebras.
It is also convenient to write the set of bicomplex numbers as
\be
   \mathbb{T}:=\{w_0+w_1\bo+w_2\bt+w_3\bj\ |\ w_0,
    w_1,w_2,w_3 \in \mathbb{R}\}.
\ee

In particular, in equation (\ref{enstetra}), if we put $z_1=x$ and
$z_2=y\bo$ with $x,y \in \mathbb{R}$, then we obtain the following
subalgebra of hyperbolic numbers, also called duplex numbers (see, e.g. \cite{14,16}):
$$\mathbb{D}:=\{x+y\bj\ |\ \bj^2=1,\ x,y\in \mathbb{R}\}\cong {\rm Cl}_{\Bbb{R}}(0,1).$$

Complex conjugation plays an important role both for algebraic and
geometric properties of $\mathbb{C}$. For bicomplex numbers, there are three possible
conjugations. Let $w\in \mathbb{T}$ and $z_1,z_2 \in
\mathbb{C}(\mathbf{i_1})$ such that $w=z_1+z_2\mathbf{i_2}$. Then we
define the three conjugations as:

\begin{subequations}
\label{eq:dag}
\begin{align}
w^{\dag_{1}}&=(z_1+z_2\bt)^{\dag_{1}}:=\overline z_1+\overline z_2
\bt,
\\
w^{\dag_{2}}&=(z_1+z_2\bt)^{\dag_{2}}:=z_1-z_2 \bt,
\\
w^{\dag_{3}}&=(z_1+z_2\bt)^{\dag_{3}}:=\overline z_1-\overline z_2
\bt,
\end{align}
\end{subequations}
where $\overline z_k$ is the standard complex conjugate of complex
numbers $z_k \in \mathbb{C}(\mathbf{i_1})$. If we say that the
bicomplex number $w=z_1+z_2\bt=w_0+w_1\bo+w_2\bt+w_3\bj$ has the
``signature'' $(++++)$, then the conjugations of type 1,2 or 3 of
$w$ have, respectively, the signatures $(+-+-)$, $(++--)$ and
$(+--+)$. We can verify easily that the composition of the
conjugates gives the four-dimensional abelian Klein group:
\begin{center}
\be
\begin{tabular}{|c||c|c|c|c|}
\hline
$\circ$ & $\dag_{0}$ & $\dag_{1}$  & $\dag_{2}$  & $\dag_{3}$  \\
\hline
\hline
$\dag_{0}$    & $\dag_{0}$ & $\dag_{1}$  & $\dag_{2}$  & $\dag_{3}$  \\
\hline
$\dag_{1}$    & $\dag_{1}$ & $\dag_{0}$ & $\dag_{3}$  & $\dag_{2}$ \\
\hline
$\dag_{2}$    & $\dag_{2}$ & $\dag_{3}$  & $\dag_{0}$ & $\dag_{1}$ \\
\hline
$\dag_{3}$    & $\dag_{3}$ & $\dag_{2}$ & $\dag_{1}$ & $\dag_{0}$ \\
\hline
\end{tabular}
\label{eq:groupedag} \ee
\end{center}
where $w^{\dag_{0}}:=w\mbox{ } \forall w\in \mathbb{T}$.

The three kinds of conjugation all have some of the standard properties of
conjugations, such as:
\begin{eqnarray}
(s+ t)^{\dag_{k}}&=&s^{\dag_{k}}+ t^{\dag_{k}},\\
\left(s^{\dag_{k}}\right)^{\dag_{k}}&=&s, \\
\left(s\cdot t\right)^{\dag_{k}}&=&s^{\dag_{k}}\cdot t^{\dag_{k}},
\end{eqnarray}
for $s,t \in \mathbb{T}$ and $k=0,1,2,3$.\\

We know that the product of a standard complex number with its
conjugate gives the square of the Euclidean metric in
$\mathbb{R}^2$. The analogs of this, for bicomplex numbers, are
the following. Let $z_1,z_2 \in \mathbb{C}(\bo)$ and
$w=z_1+z_2\bt\in \mathbb{T}$, then we have that \cite{14}:
\begin{subequations}
\begin{align}
|w|^{2}_{\bo}&:=w\cdot w^{\dag_{2}}=z^{2}_{1}+z^{2}_{2} \in
\mathbb{C}(\bo),
\\*[2ex] |w|^{2}_{\bt}&:=w\cdot
w^{\dag_{1}}=\left(|z_1|^2-|z_2|^2\right)+2\pre(z_1\overline
z_2)\bt \in \mathbb{C}(\bt), \\*[2ex] |w|^{2}_{\bj}&:=w\cdot
w^{\dag_{3}}=\left(|z_1|^2+|z_2|^2\right)-2\pim(z_1\overline
z_2)\bj \in \mathbb{D},
\end{align}
\end{subequations}
where the subscript of the square modulus refers to the subalgebra
$\mathbb{C}(\bo), \mathbb{C}(\bt)$ or $\mathbb{D}$ of $\mathbb{T}$
in which $w$ is projected.

Note that for $z_1,z_2 \in \mathbb{C}(\bo)$ and $w=z_1+z_2\bt\in
\mathbb{T}$, we can define the usual (Euclidean in $\mR^4$) norm
of $w$ as $|w|=\sqrt{|z_1|^2+|z_2|^2}=\sqrt{\pre(|w|^{2}_{\bj})}$.

It is easy to verify that $w\cdot \displaystyle
\frac{w^{\dag_{2}}}{|w|^{2}_{\bo}}=1$. Hence, the inverse of $w$
is given by \be w^{-1}= \displaystyle
\frac{w^{\dag_{2}}}{|w|^{2}_{\bo}}. \ee From this, we find that
the set $\mathcal{NC}$ of zero divisors of $\mathbb{T}$, called
the {\em null-cone}, is given by $\{z_1+z_2\bt\ |\
z_{1}^{2}+z_{2}^{2}=0\}$, which can be rewritten as \be
\mathcal{NC}=\{z(\bo\pm\bt)|\ z\in \mathbb{C}(\bo)\}. \ee
\smallskip

It is also possible to define differentiability of a function at a
point of $\mathbb{T}$:

\begin{definition}
Let $U$ be an open set of $\mathbb{T}$ and $w_0\in U$. Then,
$f:U\subseteq\mathbb{T}\longrightarrow\mathbb{T}$ is said to be
$\mathbb{T}$-differentiable at $w_{0}$ with derivative equal to
$f^\prime(w_0)\in\mathbb{T}$ if
$$\lim_{\stackrel{\scriptstyle w \rightarrow w_{0}}
{\scriptscriptstyle (w-w_{0}\mbox{
}inv.)}}\frac{f(w)-f(w_{0})}{w-w_{0}} =f^\prime(w_0).$$
\end{definition}

We also say that the function $f$ is $\mathbb{T}$-holomorphic on
an open set U if and only if $f$ is $\mathbb{T}$-differentiable at
each point of U.

\smallskip
As we saw, a bicomplex number can be seen as an element of
$\mathbb{C}^{2}$, so a function
$f(z_1+z_2\bold{i_2})=f_1(z_1,z_2)+f_2(z_1,z_2)\bold{i_2}$ of
$\mathbb{T}$ can be seen as a mapping
$f(z_1,z_2)=(f_1(z_1,z_2),f_2(z_1,z_2))$ of $\mathbb{C}^{2}$. Here
we have a characterization of such mappings:

\begin{theorem}
Let $U$ be an open set and
$f:U\subseteq\mathbb{T}\longrightarrow\mathbb{T}$ such that $f\in
{C}^{1}(U)$. Let also
$f(z_1+z_2\bold{i_2})=f_1(z_1,z_2)+f_2(z_1,z_2)\bold{i_2}$. Then
$f$ is $\mathbb{T}$-holomorphic on $U$ if and only if:
$$\mbox{$f_1$ and $f_2$ are holomorphic in $z_1$ and $z_2$}$$
and,
$$\frac{\partial{f_1}}{\partial{z_1}}
=\frac{\partial{f_2}}{\partial{z_2}}\mbox{ and
}\frac{\partial{f_2}}{\partial{z_1}}
=-\frac{\partial{f_1}}{\partial{z_2}}\mbox{ on U}.$$
\\
Moreover, $f^\prime=\frac{\partial{f_1}}{\partial{z_1}}
+\frac{\partial{f_2}}{\partial{z_1}}\bold{i_2}$ and $f^\prime(w)$ is
invertible if and only if $det\mathcal{J}_{f}(w)\neq 0$.
\label{theobasic}
\end{theorem}

\smallskip
This theorem can be obtained from results in \cite{11} and
\cite{22}. Moreover, by the Hartogs theorem \cite{21}, it is
possible to show that ``$f\in C^{1}(U)$" can be dropped from the
hypotheses. Hence, it is natural to define the
corresponding class of mappings for $\mathbb{C}^{2}$:

\begin{definition}
The class of $\mathbb{T}$-holomorphic mappings on a open set
$U\subseteq\mathbb{C}^{2}$ is defined as follows:
$$
\mbox{$TH(U):=$}\{f\mbox{:$U$}\subseteq\mathbb{C}^{2}\longrightarrow\mathbb{C}^{2}
|f\in\mbox{H($U$) and }\frac{\partial{f}_1}{\partial{z}_{1}}=
\frac{\partial{f}_2}{\partial{z}_{2}}\mbox{,
}\frac{\partial{f}_2}{\partial{z}_{1}}=-
\frac{\partial{f}_1}{\partial{z}_{2}}\mbox{ on $U$}\}.
$$

It is the subclass of holomorphic mappings of $\mathbb{C}^{2}$
satisfying the complexified Cauchy-Riemann equations.
\end{definition}

We remark that $f\in TH(U)$ in terms of $\mathbb{C}^{2}$ if and only if
$f$ is $\mathbb{T}$-differentiable on $U$.
It is also important to know that
every bicomplex number $z_1+z_2\bold{i_2}$ has the following
unique idempotent representation:
\begin{equation}
z_1+z_2\bt=(z_1-z_2\bo)\eo+(z_1+z_2\bo)\et.
\label{idempotent}
\end{equation}
where $\bold{e_1}=\frac{1+\bold{j}}{2}$ and $\bold{e_2}=\frac{1-\bold{j}}{2}$.

This representation is very useful because: addition,
multiplication and division can be done term-by-term. Also, an
element will be non-invertible if and only if $z_1-z_2\bold{i_1}=0$
or $z_1+z_2\bold{i_1}=0$.

The notion of holomorphicity can also be seen with this kind of
notation. For this we need to define the projections
$P_1,P_2:\mathbb{T}\longrightarrow\mathbb{C}(\bold{i_1})$ as
$P_1(z_1+z_2\bold{i_2})=z_1-z_2\bold{i_1}$ and
$P_2(z_1+z_2\bold{i_2})=z_1+z_2\bold{i_1}$. Also, we need the
following definition:

\begin{definition}
We say that $X\subseteq\mathbb{T}$ is a $\mathbb{T}$-cartesian set
determined by $X_1$ and $X_2$ if $X=X_{1}\times_e
X_{2}:=\{z_1+z_2\bold{i_2}\in\mathbb{T}:z_1+z_2\bold{i_2}=w_1\bold{e_1}+w_2\bold{e_2},
(w_1,w_2)\in X_1\times X_2\}$.
\end{definition}

In \cite{1} it is shown that if $X_1$ and $X_2$ are domains of
$\mathbb{C}(\bold{i_1})$ then $X_1\times_e X_2$ is also a domain
of $\mathbb{T}$. Now, it is possible to state the following
striking theorems \cite{11}:

\begin{theorem}
If $f_{e1}:X_1\longrightarrow \mathbb{C}(\bold{i_1})$ and
$f_{e2}:X_2\longrightarrow \mathbb{C}(\bold{i_1})$ are holomorphic
functions of $\mathbb{C}(\bold{i_1})$ on the domains $X_1$ and
$X_2$ respectively, then the function $f:X_1\times_e
X_2\longrightarrow \mathbb{T}$ defined as
$$f(z_1+z_2\bold{i_2})=f_{e1}(z_1-z_2\bold{i_1})\bold{e_1}+f_{e2}(z_1+z_2\bold{i_1})\bold{e_2} \mbox{ }\forall\mbox{ }z_1+z_2\bold{i_2}\in X_1\times_e X_2$$
is $\mathbb{T}$-holomorphic on the domain $X_1\times_e X_2$ and
$$f^\prime(z_1+z_2\bold{i_2})=f^\prime_{e1}(z_1-z_2\bold{i_1})\bold{e_1}+f^\prime_{e2}(z_1+z_2\bold{i_1})\bold{e_2}$$
$\mbox{ }\forall\mbox{ }z_1+z_2\bold{i_2}\in X_1\times_e X_2.$
\label{theo5}
\end{theorem}

\begin{theorem}
Let $X$ be a domain in $\mathbb{T}$, and let
$f:X\longrightarrow\mathbb{T}$ be a $\mathbb{T}$-holomorphic
function on X. Then there exist holomorphic functions
$f_{e1}:X_1\longrightarrow\mathbb{C}(\bold{i_1})$ and
$f_{e2}:X_2\longrightarrow\mathbb{C}(\bold{i_1})$ with
$X_1=P_1(X)$ and $X_2=P_2(X)$, such that:
$$f(z_1+z_2\bold{i_2})=f_{e1}(z_1-z_2\bold{i_1})\bold{e_1}+f_{e2}(z_1+z_2\bold{i_1})\bold{e_2} \mbox{ }\forall
\mbox{ }z_1+z_2\bold{i_2}\in X.$$ We note here that $X_1$ and
$X_2$ will also be domains of $\mathbb{C}(\bold{i_1})$.
\label{theo4}
\end{theorem}

\section{Bicomplex Pseudoanalytic Functions}

\subsection{Elementary Bicomplex Derivative}

We will first consider the variable $z=x+y\bo$, where $x$ and $y$
are real variables and the corresponding formal differential operators
$$\partial _{\bar{z}}  = \frac{1}{2}\left( {\partial _x  + \bo\partial _y } \right)\mbox{ and }
\partial _{z }  = \frac{1}{2}\left( {\partial _x  - \bo\partial _y } \right).$$
Notation $f_{\bar{z}}$ or $f_{z}$ means the application of $\partial _{\bar{z}}$ or $\partial _{z }$
respectively to a bicomplex function $f(z)=u(z)+v(z)\bo+r(z)\bt+s(z)\bj$. The derivatives $f_{z}$, $f_{\bar{z}}$ ``exist" if and only if
$f_{x}$ and $f_{y}$ do. Note that
\begin{eqnarray*}
f_{z} &=& \frac{1}{2}\{(u_x+v_y)+(v_x-u_y)\bo+(r_x+s_y)\bt+(s_x-r_y)\bj)\}
\end{eqnarray*}
and
\begin{eqnarray*}
f_{\bar{z}} &=& \frac{1}{2}\{(u_x-v_y)+(v_x+u_y)\bo+(r_x-s_y)\bt+(s_x+r_y)\bj)\}
\end{eqnarray*}

\noindent In view of these operators,
\begin{equation}
f_{\bar{z}}(z)=0\Leftrightarrow \partial_{\bar{z}}[u(z)+v(z)\bo]=0\mbox{ and }\partial_{\bar{z}}[r(z)+s(z)\bo]=0.
\end{equation}
$\mbox{i.e. } u_x=v_y, v_x=-u_y \mbox{ and } r_x=s_y, s_x=-r_y \mbox{ at } z\in\mathbb{C}(\bo).$

\smallskip\smallskip
We will now consider the bicomplex variable $\omega=z_1+z_2\bt$, where $z_1=x_1+y_1\bo,z_2=x_2+y_2\bo \in \mathbb{C}(\bo)$
and the corresponding formal differential operators
$$\partial _{\bar{\omega}} =\partial _{{\omega}^{\dag_{2}}}  = \frac{1}{2}\left( {\partial _{z_1}  + \bt\partial _{z_2} } \right)\mbox{, }\hspace{0.2cm}
\partial _{\omega } =\partial _{{\omega}^{\dag_{0}}} = \frac{1}{2}\left( {\partial _{z_1}  - \bt\partial _{z_2} } \right)$$

$$\partial _{{\omega}^{\dag_{3}}}  = \frac{1}{2}\left( {\partial _{\bar{z}_1}  + \bt\partial _{\bar{z}_2} } \right)\quad\mbox{ and }
\quad\partial _{{\omega}^{\dag_{1}}}  = \frac{1}{2}\left( {\partial _{\bar{z}_1}  - \bt\partial _{\bar{z}_2} } \right).$$

Notation $f_{{\omega}^{\dag_{k}}}$ for $k=0,1,2,3$ means the application of $f_{{\omega}^{\dag_{k}}}$
respectively to a bicomplex function $f(\omega)=u(\omega)+v(\omega)\bo+r(\omega)\bt+s(\omega)\bj$. The derivatives $f_{{\omega}^{\dag_{k}}}$ ``exist" for $k=0,1,2,3$ if and only if
$f_{x_l}$ and $f_{y_l}$ exist for $l=1,2$.
These bicomplex operators act on sums, products, etc. just as an ordinary derivative
and we have the following result in the bicomplex function theory.

\begin{lemma}
Let $f(z_1+z_2\bold{i_2})=f_1(z_1,z_2)+f_2(z_1,z_2)\bold{i_2}=u(z_1,z_2)+v(z_1,z_2)\bo+r(z_1,z_2)\bt+s(z_1,z_2)\bj$ be a bicomplex function. If the derivative
\begin{equation}
f'(\omega_0 ) =\lim_{\stackrel{\scriptstyle \omega \rightarrow \omega_{0}}
{\scriptscriptstyle (\omega-\omega_{0}\mbox{
}inv.)}}\frac{f(\omega)-f(\omega_{0})}{\omega-\omega_{0}}
\end{equation}
exists, then $u_x,u_y,r_x,r_y,v_x,v_y,s_x$ and $s_y$ exist, and
\begin{equation}
\mbox{1. }f_{\omega}(\omega_0)=f'(\omega_0 )\hspace{0.3cm}
\label{deriv2}
\end{equation}
\begin{equation}
\mbox{2. } f_{{\omega}^{\dag_{1}}}(\omega_0)=0\hspace{0.9cm}
\label{CR2}
\end{equation}
\begin{equation}
\mbox{3. } f_{{\omega}^{\dag_{2}}}(\omega_0)=0\hspace{1cm}
\label{CR3}
\end{equation}
\begin{equation}
\mbox{4. } f_{{\omega}^{\dag_{3}}}(\omega_0)=0.\hspace{0.9cm}
\label{CR4}
\end{equation}

Moreover, if $u_x,u_y,v_x,v_y,r_x,r_y,s_x$ and $s_y$ exist, and are continuous in a neighborhood of $\omega_0$,
and if (\ref{CR2}),(\ref{CR3}) and (\ref{CR4}) hold, then (\ref{deriv2}) exists.
\label{basic}
\end{lemma}
\emph{Proof.} First, we remark that
the complexified Cauchy-Riemann equations
\begin{equation}
\frac{\partial{f}_1}{\partial{z}_{1}}=
\frac{\partial{f}_2}{\partial{z}_{2}}\mbox{,
}\frac{\partial{f}_2}{\partial{z}_{1}}=-
\frac{\partial{f}_1}{\partial{z}_{2}}\quad \mbox{at } \omega_0
\end{equation}
are equivalent to $f_{{\omega}^{\dag_{2}}}(\omega_0)=0$. Similarly, the following system of equations
\begin{equation}
\frac{\partial{f}_1}{\partial{\bar{z}}_{1}}=
\frac{\partial{f}_2}{\partial{\bar{z}}_{2}}\mbox{,
}\frac{\partial{f}_2}{\partial{\bar{z}}_{1}}=-
\frac{\partial{f}_1}{\partial{\bar{z}}_{2}}\quad \mbox{at } \omega_0
\end{equation}
is equivalent to $f_{{\omega}^{\dag_{3}}}(\omega_0)=0$ and the following one
\begin{equation}
\frac{\partial{f}_1}{\partial{\bar{z}}_{1}}=-
\frac{\partial{f}_2}{\partial{\bar{z}}_{2}}\mbox{,
}\frac{\partial{f}_2}{\partial{\bar{z}}_{1}}=
\frac{\partial{f}_1}{\partial{\bar{z}}_{2}}\quad \mbox{at } \omega_0
\end{equation}
is equivalent to $f_{{\omega}^{\dag_{1}}}(\omega_0)=0$.
Now, if the bicomplex derivative $f'(\omega_0 )$ exists, then from the Theorem \ref{theobasic}
we obtain automatically that all partial derivatives exist at $\omega_0$.
Moreover, the fact that $\mbox{$f_1$ and $f_2$ are holomorphic in $z_1$ and $z_2$}$ imply that
\begin{equation}
\frac{\partial{f}_k}{\partial{\bar{z}}_{1}}=\frac{\partial{f}_k}{\partial{\bar{z}}_{2}}=0\quad \mbox{for }k=1,2.
\end{equation}
Therefore $f_{{\omega}^{\dag_{1}}}(\omega_0)=f_{{\omega}^{\dag_{3}}}(\omega_0)=0$ and the complexified Cauchy-Riemann equations
imply $f_{{\omega}^{\dag_{2}}}(\omega_0)=0$ and $f_{\omega}(\omega_0)=f'(\omega_0 )$. Conversely, if $u_x,u_y,r_x,r_y,v_x,v_y$, $s_x$ and $s_y$ exist,
and are continuous in a neighborhood of $\omega_0$ then $f_{{\omega}^{\dag_{1}}}(\omega_0)=f_{{\omega}^{\dag_{3}}}(\omega_0)=0$ imply that
$\mbox{$f_1$ and $f_2$ are holomorphic in $z_1$ and $z_2$}$. Hence, from Theorem \ref{theobasic}, $f$ is $\mathbb{T}$-differentiable
at $\omega_0$.$\Box$

\subsection{Bicomplex Generalization of Function Theory}

Our bicomplex generalization of function theory is based on the following three different representations of
bicomplex numbers. The \textit{scalar} and \textit{vectorial} part must be adapted to each representations.

\subsubsection{Class-R1}

Let $a+b\bo+c\bt+d\bj=z_1+z_2\bt$ where $z_1,z_2\in\mathbb{C}(\bo)$. In this case, the
theory will be based on assigning the part played by $1$ and $\bt$ to two essentially arbitrary bicomplex functions
$F(\omega)$ and $G(\omega)$. We assume that these functions are defined and twice continuously differentiable in some open domain $D_{0}\subset\mathbb{T}$. We require that
\begin{equation}
\mbox{Vec}\{F(\omega)^{\dag_{2}}G(\omega)\}\neq0.
\label{vec01}
\end{equation}
Under this condition, $(F,G)$ will be called a ${\bo}$-generating pair in $D_0$.
We remark that $\mbox{Vec}\{F(\omega)^{\dag_{2}}G(\omega)\}=
\left| {\begin{array}{*{20}c}
{\mbox{Sc}\{F(\omega)\}} & {\mbox{Sc}\{G(\omega)\}}  \\
   {\mbox{Vec}\{F(\omega)\}} & {\mbox{Vec}\{G(\omega)\}}  \\
\end{array} } \right|.$
It follows, from Cramer's Theorem, that for every $\omega_0$ in $D_0$ we can find \textbf{unique} constants $\lambda_0,\mu_0\in\mathbb{C}(\bo)$ such that
$w(\omega_0)=\lambda_0 F(\omega_0)+\mu_0 G(\omega_0)$. More generally we have the following result.

\begin{theorem}
Let $(F,G)$ be $\bo$-generating pair in some open domain $D_{0}$.
If  $w(\omega):D_{0}\subset\mathbb{T}\to\mathbb{T}$, then there exist \textbf{unique} functions
$\phi(\omega), \psi(\omega):D_{0}\subset\mathbb{T}\to\mathbb{C}(\bo)$ such that
$$w(\omega)=\phi(\omega)F(\omega)+\psi(\omega)G(\omega)\mbox{ }\forall \omega\in D_{0}.$$ Moreover, we have the following explicit
formulas for $\phi$ and  $\psi$:
$$\phi(\omega)=\frac{Vec[w(\omega)^{\dag_{2}}G(\omega)]}{Vec[F(\omega)^{\dag_{2}}G(\omega)]}\mbox{, }
\psi(\omega)=-\frac{Vec[w(\omega)^{\dag_{2}}F(\omega)]}{Vec[F(\omega)^{\dag_{2}}G(\omega)]}.$$
\label{explicit}
\end{theorem}
\emph{Proof.}
Let $(F,G)$ be $\bo$-generating pair in some open domain $D_{0}$.
Let $z_0\in D_{0}$ with $w(z_0)=z_1+z_2\bt$, $F(z_0)=z_3+z_4\bt$ and $G(z_0)=z_5+z_6\bt$.
In this case, $w(z_0)=\phi_{2}(z_0)F(z_0)+\psi_{2}(z_0)G(z_0)$ with $\phi_{2}(z_0), \psi_{2}(z_0)\in\mathbb{C}(\bo)$
if and only if $z_1=\phi_{2}(z_0)z_3+\psi_{2}(z_0)z_5$ and $z_2=\phi_{2}(z_0)z_4+\psi_{2}(z_0)z_6$. This is a well known
Cramer's system of the form $AX=B$ where $A = \left( {\begin{array}{*{20}c}
   {z_3 } & {z_5 }  \\
   {z_4 } & {z_6 }  \\
\end{array} } \right)$,
$B = \left( {\begin{array}{*{20}c}
   {z_1 }  \\
   {z_2 }  \\
\end{array} } \right)$ and
$X = \left( {\begin{array}{*{20}c}
   {\phi_{2}(z_0 )}  \\
   {\psi_{2}(z_0 )}  \\
\end{array} } \right)$.
So, the unique solution is $X=A^{-1}B$ where $A^{-1}=
\frac{1}
{{\det A}}\left( {\begin{array}{*{20}c}
   {z_6 } & { - z_5 }  \\
   { - z_4 } & {z_3 }  \\
\end{array} } \right)$. Hence, $X=
\frac{1}
{{\mbox{Vec}[F(z_0)^{\dag_{2}}G(z_0)]}}\left( {\begin{array}{*{20}c}
   {z_6 } & { - z_5 }  \\
   { - z_4 } & {z_3 }  \\
\end{array} } \right)
\left( {\begin{array}{*{20}c}
   {z_1 }  \\
   {z_2 }  \\
\end{array} } \right)=
\frac{1}
{{\mbox{Vec}[F(z_0)^{\dag_{2}}G(z_0)]}}\left( {\begin{array}{*{20}c}
   {z_1z_6 } - {z_2z_5}  \\
   {-z_1z_4} + {z_2z_3}  \\
\end{array} } \right)\\
=\frac{1}
{{\mbox{Vec}[F(z_0)^{\dag_{2}}G(z_0)]}}\left( {\begin{array}{*{20}c}
   {\mbox{Vec}[w(z_0)^{\dag_{2}}G(z_0)]}  \\
   {-\mbox{Vec}[w(z_0)^{\dag_{2}}F(z_0)]}  \\
\end{array} }\right)$.\\
Then
$$\phi_{2}(z)=\frac{\mbox{Vec}[w(z)^{\dag_{2}}G(z)]}{\mbox{Vec}[F(z)^{\dag_{2}}G(z)]}\mbox{, }
\psi_{2}(z)=-\frac{\mbox{Vec}[w(z)^{\dag_{2}}F(z)]}{\mbox{Vec}[F(z)^{\dag_{2}}G(z)]}\mbox{ }\forall z\in D_0.\Box$$

\smallskip\smallskip
\hspace{0.5cm} Now, we say that $w(\omega):D_{0}\subset\mathbb{T}\to\mathbb{T}$ possesses at $\omega_0$ the $(F,G)_{\bo}$-derivative $\dot w(\omega_0)$ if the (finite) limit
\begin{equation}
\dot w(\omega_0)=\lim_{\stackrel{\scriptstyle \omega \rightarrow \omega_{0}}{\scriptscriptstyle (\omega-\omega_{0}\mbox{
}inv.)}} \frac{{w(\omega) - \lambda _0 F(\omega) - \mu _0 G(\omega)}}
{{\omega - \omega_0 }}
\label{lim1}
\end{equation}
exists.
\smallskip

In the particular case where $w(\omega),F(\omega)$ and $G(\omega)$ are defined
on $D_{0}\subset\mathbb{C}(\bt)\to\mathbb{C}(\bt)$ then
we can find unique constants $\lambda_0,\mu_0\in\mathbb{R}$ such that
$w(\omega_0)=\lambda_0 F(\omega_0)+\mu_0 G(\omega_0)$ and we come back to the classical (in $\bt$)
pseudoanalytic theory developed by L. Bers and I.N. Vekua (see e.g. \cite{1,3,4,17}). In that case, using Bers's theory of Taylor series for pseudoanalytic function,
V.V. Kravchenko (see \cite{8}) obtained a locally complete system of solutions of the real stationary two-dimensional Schr{\"o}dinger equation.
On the other hand, in the case where
$w(\omega)$, $F(\omega)$ and $G(\omega)$ are defined on $D_{0}\subset\mathbb{C}(\bj)\to\mathbb{C}(\bj)$ then
we can also find unique constants $\lambda_0,\mu_0\in\mathbb{R}$ such that
$w(\omega_0)=\lambda_0 F(\omega_0)+\mu_0 G(\omega_0)$ and we are in the hyperbolic
pseudoanalytic theory developed by Guo Chun Wen in \cite{18}.
Moreover, we note that if we only restrict the domain $D_0$ to $\mathbb{C}(\bt)$, the subclass of bicomplex pseudoanalytic functions obtained
is precisely the class developed by V.V. Kravchenko and A. Casta\~neda in \cite{5} to show that in a two-dimensional situation the Dirac equation
with a scalar and an electromagnetic potentials decouples into a pair of bicomplex equations. It is also the same class of functions
that used V.V. Kravchenko (see \cite{20}) to obtain solutions of the complex stationary two-dimensional Schr{\"o}dinger equation.

\hspace{0.5cm} The following expressions are called the $\bo$-characteristic coefficients of the pair $(F,G)$ for $k=1,2,3$:
$$a_{(F,G)}^{(k)} =-\frac{F^{\dag_{k}}G_{\omega^{\dag_{k}}}-F_{\omega^{\dag_{k}}}G^{\dag_{k}}}{FG^{\dag_{2}}-F^{\dag_{2}}G},\hspace{1cm}
b_{(F,G)}^{(k)} =\frac{FG_{\omega^{\dag_{k}}}-F_{\omega^{\dag_{k}}}G}{FG^{\dag_{2}}-F^{\dag_{2}}G},$$

$$A_{(F,G)}=-\frac{F^{\dag_{2}}G_{\omega}-F_{\omega}G^{\dag_{2}}}{FG^{\dag_{2}}-F^{\dag_{2}}G},\hspace{1cm}
B_{(F,G)}=\frac{FG_{\omega}-F_{\omega}G}{FG^{\dag_{2}}-F^{\dag_{2}}G}.$$

Set (for a fixed $\omega_0$)
\begin{equation}
W(\omega)=w(\omega)-\lambda _0 F(\omega) - \mu _0 G(\omega),
\end{equation}
the constants $\lambda_0,\mu_0\in\mathbb{C}(\bo)$ being uniquely determined by the condition
\begin{equation}
W(\omega_0)=0.
\end{equation}
Hence, $W(\omega)$ has continuous partial derivatives if and only if $w(\omega)$ has. Moreover, $\dot w(\omega_0)$ exists if and only if $W'(\omega_0)$ does, and if
it does exist, then $\dot w(\omega_0)=W'(\omega_0)$. Therefore, by Lemma \ref{basic}, the existence of $W_\omega(\omega_0)$, $W_{\bar{\omega}}(\omega_0)$ and equations
\begin{equation}
W_{\omega^{\dag_{1}}}(\omega_0)=0,
 W_{\bar{\omega}}(\omega_0)=0 \mbox{ and }
W_{\omega^{\dag_{3}}}(\omega_0)=0
\label{cod}
\end{equation}
are necessary for the existence of (\ref{lim1}), and the existence and continuity of $W_\omega(\omega_0)$, $W_{\bar{\omega}}(\omega_0)$
in a neighborhood of $\omega_0$, together with (\ref{cod}) are sufficient. Now,
\begin{equation}
W(\omega) = \frac{{\left| {\begin{array}{*{20}c}
   {w(\omega)} & {w(\omega_0 )} & {w(\omega_0 )^{\dag_{2}}}  \\
   {F(\omega)} & {F(\omega_0 )} & {F(\omega_0 )^{\dag_{2}}}  \\
   {G(\omega)} & {G(\omega_0 )} & {G(\omega_0 )^{\dag_{2}}}  \\

 \end{array} } \right|}}
{{\left| {\begin{array}{*{20}c}
   {F(\omega_0 )} & {F(\omega_0 )^{\dag_{2}}}  \\
   {G(\omega_0 )} & {G(\omega_0 )^{\dag_{2}}}  \\

 \end{array} } \right|}}
\end{equation}

so that (\ref{cod}) may be written in the form
\begin{equation}
\left| {\begin{array}{*{20}c}
   {w_{\omega^{\dag_{k}}}(\omega_0)} & {w(\omega_0 )} & {w(\omega_0 )^{\dag_{2}}}  \\
   {F_{\omega^{\dag_{k}}}(\omega_0)} & {F(\omega_0 )} & {F(\omega_0 )^{\dag_{2}}}  \\
   {G_{\omega^{\dag_{k}}}(\omega_0)} & {G(\omega_0 )} & {G(\omega_0 )^{\dag_{2}}}  \\
 \end{array} } \right| = 0 \quad\mbox{for } k=1,2,3
\label{E1}
\end{equation}
and if (\ref{lim1}) exists, then
\begin{equation}
 \dot w(\omega_0)= \frac{{\left| {\begin{array}{*{20}c}
   {w_\omega(\omega_0)} & {w(\omega_0 )} & {w(\omega_0 )^{\dag_{2}}}  \\
   {F_\omega(\omega_0)} & {F(\omega_0 )} & {F(\omega_0 )^{\dag_{2}}}  \\
   {G_\omega(\omega_0)} & {G(\omega_0 )} & {G(\omega_0 )^{\dag_{2}}}  \\

 \end{array} } \right|}}
{{\left| {\begin{array}{*{20}c}
   {F(\omega_0 )} & {F(\omega_0 )^{\dag_{2}}}  \\
   {G(\omega_0 )} & {G(\omega_0 )^{\dag_{2}}}  \\

 \end{array} } \right|}}.
\label{E2}
\end{equation}
Equations (\ref{E2}) and (\ref{E1}) can be rewritten in the form
\begin{equation}
\dot w=w_{\omega}-A_{(F,G)}w-B_{(F,G)}w^{\dag_{2}}
\end{equation}
\begin{equation}
w_{\omega^{\dag_{k}}}=a_{(F,G)}^{(k)} w+b_{(F,G)}^{(k)} w^{\dag_{2}}\quad\mbox{for } k=1,2,3.
\end{equation}

Thus we have proved the following result.
\begin{theorem}
Let (F,G) be a ${\bo}$-generating pair in some open domain $D_0$. Every bicomplex function $w$ defined in $D_0$ admits the unique representation $w=\phi F + \psi G$ where
$\phi,\psi:D_{0}\subset\mathbb{T}\to\mathbb{C}(\bo)$. Moreover, the $(F,G)_{\bo}$-derivative
$\dot w=\frac{d_{(F,G)_{\bo}} w}{d\omega}$ of $w(\omega)$ exists at $\omega_0$ and has the form
\begin{equation}
\dot w=\phi_{\omega}F+\psi_{\omega}G=w_{\omega}-A_{(F,G)}w-B_{(F,G)}w^{\dag_{2}}
\end{equation}
if and only if
\begin{equation}
w_{\omega^{\dag_{1}}}=a_{(F,G)}^{(1)}w+b_{(F,G)}^{(1)}w^{\dag_{2}},
\label{vekua11}
\end{equation}
\begin{equation}
w_{\omega^{\dag_{2}}}=a_{(F,G)}^{(2)}w+b_{(F,G)}^{(2)}w^{\dag_{2}},
\label{vekua21}
\end{equation}
and
\begin{equation}
w_{\omega^{\dag_{3}}}=a_{(F,G)}^{(3)}w+b_{(F,G)}^{(3)}w^{\dag_{2}}
\label{vekua31}
\end{equation}
where $w$ has continuous partial derivatives in a neighborhood of $\omega_0$.
\label{iifR1}
\end{theorem}

\smallskip
The equations (\ref{vekua11}), (\ref{vekua21}) and (\ref{vekua31}) are called the $\bo$-bicomplex Vekua equations and the solutions of
these equations will be the $(F,G)_{\bo}$-pseudoanalytic functions.

\begin{remark}
For $k=1,2,3$, the equation
\begin{equation}
w_{\omega^{\dag_{k}}}=a_{(F,G)}^{(k)}w+b_{(F,G)}^{(k)}w^{\dag_{2}}
\label{equi01}
\end{equation}
can be rewritten in the following form
\begin{equation}
\phi_{\omega^{\dag_{k}}}F+\psi_{\omega^{\dag_{k}}}G=a_{(F,G)}^{(k)}w +\frac{F^{\dag_{2}}G_{\omega^{\dag_{k}}}-F_{\omega^{\dag_{k}}}G^{\dag_{2}}}{FG^{\dag_{2}}-F^{\dag_{2}}G}w.
\end{equation}
Hence, the equation (\ref{equi01}) is equivalent to $$\phi_{\omega^{\dag_{k}}}F+\psi_{\omega^{\dag_{k}}}G=0$$ if and only if
\begin{equation}
[G^{\dag_{k}}-G^{\dag_{2}}]F_{\omega^{\dag_{k}}}=[F^{\dag_{k}}-F^{\dag_{2}}]G_{\omega^{\dag_{k}}}
\end{equation}
where $w$ is not identically in the null-cone on the domain.
\end{remark}

\subsubsection{Class-R2}

Let $a+b\bo+c\bt+d\bj=z_1+z_2\bo$ where $z_1,z_2\in\mathbb{C}(\bt)$. In this case, the
theory will be based on assigning the part played by $1$ and $\bo$ to two essentially arbitrary bicomplex functions
$F(\omega)$ and $G(\omega)$. We assume that these functions are defined and twice continuously differentiable in some open domain $D_{0}\subset\mathbb{T}$. We require that
\begin{equation}
\mbox{Vec}\{F(\omega)^{\dag_{1}}G(\omega)\}\neq0.
\label{vec02}
\end{equation}
Under this condition, $(F,G)$ will be called a ${\bt}$-generating pair in $D_0$.
From Cramer's Theorem, it follows that for every $\omega_0$ in $D_0$ we can find \textbf{unique} constants $\lambda_0,\mu_0\in\mathbb{C}(\bt)$ such that
$w(\omega_0)=\lambda_0 F(\omega_0)+\mu_0 G(\omega_0)$. In fact, using the same arguments than for the Theorem \ref{explicit}, we have the following result.

\begin{theorem}
Let $(F,G)$ be $\bt$-generating pair in some open domain $D_{0}$.
If  $w(\omega):D_{0}\subset\mathbb{T}\to\mathbb{T}$, then there exist \textbf{unique} functions
$\phi(\omega), \psi(\omega):D_{0}\subset\mathbb{T}\to\mathbb{C}(\bt)$ such that
$$w(\omega)=\phi(\omega)F(\omega)+\psi(\omega)G(\omega)\mbox{ }\forall \omega\in D_{0}.$$ Moreover, we have the following explicit
formulas for $\phi$ and  $\psi$:
$$\phi(\omega)=\frac{Vec[w(\omega)^{\dag_{1}}G(\omega)]}{Vec[F(\omega)^{\dag_{1}}G(\omega)]}\mbox{, }
\psi(\omega)=-\frac{Vec[w(\omega)^{\dag_{1}}F(\omega)]}{Vec[F(\omega)^{\dag_{1}}G(\omega)]}.$$
\end{theorem}

\smallskip\smallskip
\hspace{0.5cm} We say that $w(\omega):D_{0}\subset\mathbb{T}\to\mathbb{T}$ possesses at $\omega_0$ the $(F,G)_{\bt}$-derivative $\dot w(\omega_0)$
if the (finite) limit
\begin{equation}
\dot w(\omega_0)=\lim_{\stackrel{\scriptstyle \omega \rightarrow \omega_{0}}{\scriptscriptstyle (\omega-\omega_{0}\mbox{
}inv.)}} \frac{{w(\omega) - \lambda _0 F(\omega) - \mu _0 G(\omega)}}
{{\omega - \omega_0 }}
\label{lim2}
\end{equation}
exists. In fact, if we interchange everywhere $\bo$ with $\bt$, this case is exactly the same than \textbf{R1}.
In particular, if we defined the function $\pi:\mathbb{T}\longrightarrow\mathbb{T}$ as
\begin{equation}
\pi(a+b\bo+c\bt+d\bj):=a+c\bo+b\bt+d\bj
\end{equation}
we obtain that $w(\omega)$ possesses a $(F,G)_{\bo}$-derivative at $\omega_0\in D_0$ if and only if
the function
\begin{equation}
(\pi\circ w\circ\pi)(\omega)
\end{equation}
possesses a $(\pi\circ F\circ\pi,\pi\circ G\circ\pi)_{\bt}$-derivative at $\pi(\omega_0)\in\pi(D_0)$
where
\begin{equation}
(\pi\circ F\circ\pi,\pi\circ G\circ\pi)
\end{equation}
is a ${\bt}$-generating pair on $\pi(D_0)$. We note that
\begin{eqnarray}
\pi(\pi(\omega)) &=& \pi(\omega),\\
\pi(\omega_1 + \omega_2) &=& \pi(\omega_1)+\pi(\omega_2),\\
\pi(\omega_1 \omega_2) &=& \pi(\omega_1)\pi(\omega_2),\\
\pi\left(\frac{\omega_1}{\omega_2}\right) &=& \frac{\pi(\omega_1)}{\pi(\omega_2)}\quad \mbox{if }\omega_2\notin\mathcal{NC},\\
\pi(w^{\dag_{3}}) &=& (\pi(w))^{\dag_{3}},\\
\pi(w^{\dag_{1}}) = (\pi(w))^{\dag_{2}} &\mbox{ and }& \pi(w^{\dag_{2}}) = (\pi(w))^{\dag_{1}}.
\end{eqnarray}
Therefore,
$\left.{\frac{d_{(F,G)_{\bt}} {(\pi\circ w\circ\pi)}}{d\omega}}\right|_{\omega=\pi(\omega_0)}$
\small
\begin{eqnarray*}
 &=& \lim_{\stackrel{\scriptstyle \omega \rightarrow \pi(\omega_{0})}{\scriptscriptstyle
(\omega-\pi(\omega_{0})\mbox{ }inv.)}} \frac{{(\pi\circ w\circ\pi)(\omega) - \pi(\lambda _0) (\pi\circ F\circ\pi)(\omega) - \pi(\mu _0) (\pi\circ G\circ\pi)(\omega)}}{{\omega - \pi(\omega_0 })}\\
&=& \lim_{\stackrel{\scriptstyle \omega \rightarrow \pi(\omega_{0})}{\scriptscriptstyle
(\omega-\pi(\omega_{0})\mbox{ }inv.)}} \frac{\pi[w(\pi(\omega)) - \lambda _0 F(\pi(\omega)) - \mu _0 G(\pi(\omega))]}{\pi[{\pi(\omega) - \omega_0}]}\\
&=& \pi \left[ \lim_{\stackrel{\scriptstyle \omega \rightarrow \pi(\omega_{0})}{\scriptscriptstyle
(\omega-\pi(\omega_{0})\mbox{ }inv.)}} \frac{w(\pi(\omega)) - \lambda _0 F(\pi(\omega)) - \mu _0 G(\pi(\omega))}{{\pi(\omega) - \omega_0}} \right]\\
&=& \pi \left[ \left.{\frac{d_{(F,G)_{\bo}} w }{d\omega}}\right|_{\omega=\omega_0} \right].
\end{eqnarray*}
\normalsize
\hspace{0.5cm} We note that the $\bt$-characteristic coefficients of the pair $(F,G)$ for $k=1,2,3$ must be defined as:
$$a_{(F,G)}^{(k)} =-\frac{F^{\dag_{k}}G_{\omega^{\dag_{k}}}-F_{\omega^{\dag_{k}}}G^{\dag_{k}}}{FG^{\dag_{1}}-F^{\dag_{1}}G},\hspace{1cm}
b_{(F,G)}^{(k)} =\frac{FG_{\omega^{\dag_{k}}}-F_{\omega^{\dag_{k}}}G}{FG^{\dag_{1}}-F^{\dag_{1}}G},$$

$$A_{(F,G)}=-\frac{F^{\dag_{1}}G_{\omega}-F_{\omega}G^{\dag_{1}}}{FG^{\dag_{1}}-F^{\dag_{1}}G},\hspace{1cm}
B_{(F,G)}=\frac{FG_{\omega}-F_{\omega}G}{FG^{\dag_{1}}-F^{\dag_{1}}G}.$$

\subsubsection{Class-R3}

Let $a+b\bo+c\bt+d\bj=z_1+z_2\bo$ (resp. $z_1+z_3\bt$) where $z_1,z_2,z_3\in\mathbb{C}(\bj)$. In this case, the
theory will be based on assigning the part played by $1$ and $\bo$ (resp. $\bt$) to two essentially arbitrary bicomplex functions
$F(z)$ and $G(z)$. We assume that these functions are defined and twice continuously differentiable in some open domain $D_{0}\subset\mathbb{T}$. We require that
\begin{equation}
\mbox{Vec}\{F(\omega)^{\dag_{3}}G(\omega)\}\neq0.
\label{vec03}
\end{equation}
Under this condition, $(F,G)$ will be called a ${\bj}$-generating pair in $D_0$
for every $\omega_0$ in $D_0$. Moreover, it will be possible to find \textbf{unique} constants $\lambda_0,\mu_0\in\mathbb{C}(\bj)$ such that
$w(\omega_0)=\lambda_0 F(\omega_0)+\mu_0 G(\omega_0)$. Here also we have the following equivalence of the Theorem \ref{explicit}.

\begin{theorem}
Let $(F,G)$ be $\bj$-generating pair in some open domain $D_{0}$.
If  $w(\omega):D_{0}\subset\mathbb{T}\to\mathbb{T}$, then there exist \textbf{unique} functions
$\phi(\omega), \psi(\omega):D_{0}\subset\mathbb{T}\to\mathbb{C}(\bj)$ such that
$$w(\omega)=\phi(\omega)F(\omega)+\psi(\omega)G(\omega)\mbox{ }\forall \omega\in D_{0}.$$ Moreover, we have the following explicit
formulas for $\phi$ and  $\psi$:
$$\phi(\omega)=\frac{Vec[w(\omega)^{\dag_{3}}G(\omega)]}{Vec[F(\omega)^{\dag_{3}}G(\omega)]}\mbox{, }
\psi(\omega)=-\frac{Vec[w(\omega)^{\dag_{3}}F(\omega)]}{Vec[F(\omega)^{\dag_{3}}G(\omega)]}.$$
\end{theorem}

\smallskip\smallskip
\hspace{0.5cm} We say that $w(\omega)$ possesses at $\omega_0$ the $(F,G)_{\bj}$-derivative $\dot w(\omega_0)$ if the (finite) limit
\begin{equation}
\dot w(\omega_0)=\mathop {\lim }\limits_{\omega \to \omega_0 } \frac{{w(\omega) - \lambda _0 F(\omega) - \mu _0 G(\omega)}}
{{\omega - \omega_0 }}
\end{equation}
exists.

\hspace{0.5cm} In this case, the following expressions are called the $\bj$-characteristic coefficients of the pair $(F,G)$ for $k=1,2,3$:
$$a_{(F,G)}^{(k)} =-\frac{F^{\dag_{k}}G_{\omega^{\dag_{k}}}-F_{\omega^{\dag_{k}}}G^{\dag_{k}}}{FG^{\dag_{3}}-F^{\dag_{3}}G},\hspace{1cm}
b_{(F,G)}^{(k)} =\frac{FG_{\omega^{\dag_{k}}}-F_{\omega^{\dag_{k}}}G}{FG^{\dag_{3}}-F^{\dag_{3}}G},$$

$$A_{(F,G)}=-\frac{F^{\dag_{3}}G_{\omega}-F_{\omega}G^{\dag_{3}}}{FG^{\dag_{3}}-F^{\dag_{3}}G},\hspace{1cm}
B_{(F,G)}=\frac{FG_{\omega}-F_{\omega}G}{FG^{\dag_{3}}-F^{\dag_{3}}G}.$$

Now, using the same kind of arguments than for the case \textbf{R2}, we obtain the following result.
\begin{theorem}
Let (F,G) be a ${\bj}$-generating pair in some open domain $D_0$. Every bicomplex function $w$ defined in $D_0$ admits the unique representation $w=\phi F + \psi G$ where
$\phi,\psi:D_{0}\subset\mathbb{T}\to\mathbb{C}(\bj)$. Moreover, the $(F,G)_{\bj}$-derivative
$\dot w=\frac{d_{(F,G)_{\bj}} w}{d\omega}$ of $w(\omega)$ exists at $\omega_0$ and has the form
\begin{equation}
\dot w=\phi_{\omega}F+\psi_{\omega}G=w_{\omega}-A_{(F,G)}w-B_{(F,G)}w^{\dag_{3}}
\end{equation}
if and only if
\begin{equation}
w_{\omega^{\dag_{1}}}=a_{(F,G)}^{(1)}w+b_{(F,G)}^{(1)}w^{\dag_{3}},
\label{vekua13}
\end{equation}
\begin{equation}
w_{\omega^{\dag_{2}}}=a_{(F,G)}^{(2)}w+b_{(F,G)}^{(2)}w^{\dag_{3}},
\label{vekua23}
\end{equation}
and
\begin{equation}
w_{\omega^{\dag_{3}}}=a_{(F,G)}^{(3)}w+b_{(F,G)}^{(3)}w^{\dag_{3}}
\label{vekua33}
\end{equation}
where $w$ has continuous partial derivatives in a neighborhood of $\omega_0$.
\label{iifR3}
\end{theorem}

\smallskip
The equations (\ref{vekua13}), (\ref{vekua23}) and (\ref{vekua33}) are called the $\bj$-bicomplex Vekua equations and the solutions of
these equations will be the $(F,G)_{\bj}$-pseudoanalytic functions.

\begin{remark}
For $k=1,2,3$, the equation
\begin{equation}
w_{\omega^{\dag_{k}}}=a_{(F,G)}^{(k)}w+b_{(F,G)}^{(k)}w^{\dag_{3}}
\label{equi03}
\end{equation}
can be rewritten in the following form
\begin{equation}
\phi_{\omega^{\dag_{k}}}F+\psi_{\omega^{\dag_{k}}}G=a_{(F,G)}^{(k)}w +\frac{F^{\dag_{3}}G_{\omega^{\dag_{k}}}-F_{\omega^{\dag_{k}}}G^{\dag_{3}}}{FG^{\dag_{3}}-F^{\dag_{3}}G}w.
\end{equation}
Hence, the equation (\ref{equi03}) is equivalent to $$\phi_{\omega^{\dag_{k}}}F+\psi_{\omega^{\dag_{k}}}G=0$$ if and only if
\begin{equation}
[G^{\dag_{k}}-G^{\dag_{3}}]F_{\omega^{\dag_{k}}}=[F^{\dag_{k}}-F^{\dag_{3}}]G_{\omega^{\dag_{k}}}
\end{equation}
where $w$ is not identically in the null-cone on the domain.
\end{remark}

\hspace{0.5cm} In this case, it is useful to consider a more specific class of generating pair.
\begin{definition}
Let $D_1$ and $D_2$ be open in $\mathbb{C}(\bo)$. Consider that $(F_{e_1},G_{e_1})$ and $(F_{e_2},G_{e_2})$ are complex (in $\bo$),
twice continuously differentiable, generating pairs
in respectively $D_1$ and $D_2$. Under these conditions, $(F,G)$ will be called a
${\bje}$-generating pair in $D_0=D_1\times_e D_2\in\mathbb{T}$ where
\begin{equation}
F(z_1+z_2\bt):=F_{e_1}(z_1-z_2\bo)\eo+F_{e_2}(z_1+z_2\bo)\et
\end{equation}
and
\begin{equation}
G(z_1+z_2\bt):=G_{e_1}(z_1-z_2\bo)\eo+G_{e_2}(z_1+z_2\bo)\et.
\end{equation}
\end{definition}

\begin{lemma}
Let $F(\omega)$ and $G(\omega)$ two arbitrary bicomplex functions defined in some domain $D_{0}\subset\mathbb{T}$.
If $$\mbox{Im}\{\overline{F_{e_1}(\omega)}G_{e_1}(\omega)\}\neq0\mbox{ or }\mbox{Im}\{ \overline{F_{e_2}(\omega)}G_{e_2}(\omega)\}\neq0
\mbox{ }\forall \omega\in D_0$$ then
$\mbox{Vec}\{F(\omega)^{\dag_{3}}G(\omega)\}\neq0\mbox{ }\forall \omega\in D_0$.
\label{lem01}
\end{lemma}
\emph{Proof.}
Let $F(\omega_0)^{\dag_{3}}G(\omega_0)=x_1+x_2\bo+x_3\bt+x_4\bj=(x_1+x_4\bj)+(x_2-x_3\bj)\bo=(x_1+x_4\bj)+(x_3-x_2\bj)\bt$. Thus,
$\mbox{Vec}\{F(\omega_0)^{\dag_{3}}G(\omega_0)\}=0$ if and only if $x_2=x_3=0$. Moreover,
$\mbox{Im}\{ \overline{F_{e_1}(\omega_0)}G_{e_1}(z_0)\}=x_3-x_2$ and $\mbox{Im}\{ \overline{F_{e_2}(\omega_0)}G_{e_2}(\omega_0)\}=x_2+x_3$.
Hence, $\mbox{Vec}\{F(\omega_0)^{\dag_{3}}G(\omega_0)\}=0$ imply $\mbox{Im}\{ \overline{F_{e_1}(z_0)}G_{e_1}(\omega_0)\}=0$ and
$\mbox{Im}\{ \overline{F_{e_2}(\omega_0)}G_{e_2}(\omega_0)\}=0\mbox{ }\forall \omega_0\in D_0$.$\Box$
\smallskip\smallskip

\hspace{0.5cm} Therefore, from the Lemma \ref{lem01} we obtain automatically this following result.
\begin{theorem}
Let $D_0=D_1\times_e D_2$ where $D_1$ and $D_2$ are open domains in $\mathbb{C}(\bo)$.
If $(F,G)$ is a ${\bje}$-generating pair in $D_0$ then $(F,G)$ is, in particular,
a ${\bj}$-generating pair in $D_0$.
\end{theorem}

\hspace{0.5cm} Moreover, the ${\bje}$-generating pair will imply the following representation for a
$(F,G)_{\bj}$-pseudoanalytic function.
\begin{theorem}
Let $D_0=D_1\times_e D_2$ where $D_1$ and $D_2$ are open domains in $\mathbb{C}(\bo)$.
If $(F,G)$ is a ${\bje}$-generating pair in $D_0$ then if the function $w$ is $(F,G)_{\bj}$-pseudoanalytic
on $D_0$ then $w$ can be decomposed in the following way
\begin{equation}
w(z_1+z_2\bt)=[{w}_{e_1}(z_1-z_2\bo)]\eo+[{w}_{e_2}(z_1+z_2\bo)]\et
\end{equation}
$\forall$ $z_1+z_2\bt=(z_1-z_2\bo)\eo+(z_1+z_2\bo)\et\in D_0$.
\end{theorem}
\emph{Proof.}
It is always possible to decomposed $w(z_1+z_2\bt)$ in term of the idempotent representation i.e.
$$w(z_1+z_2\bt)=[{w}_{e_1}(z_1+z_2\bt)]\eo+[{w}_{e_2}(z_1+z_2\bt)]\et.$$
Moreover, from the definition of the derivative, the function
$$W(z_1+z_2\bt)=w(z_1+z_2\bt)-\lambda _0 F(z_1+z_2\bt) - \mu _0 G(z_1+z_2\bt)$$
is $\mathbb{T}$-differentiable at $z_1+z_2\bt$. Now, using the Theorem \ref{theo5},
we have that
$$W(z_1+z_2\bt)=[{W}_{e_1}(z_1-z_2\bo)]\eo+[{W}_{e_2}(z_1+z_2\bo)]\et$$
and from the definition of a ${\bje}$-generating pair we obtain that
\begin{eqnarray*}
w(z_1+z_2\bt) &=& {W}_{e_1}(z_1-z_2\bo)\eo+{W}_{e_2}(z_1+z_2\bo)\et\\
              &+& P_1(\lambda_0)F_{e_1}(z_1-z_2\bo)\eo+P_2(\lambda_0)F_{e_2}(z_1+z_2\bo)\et\\
              &+& P_1(\mu_0)G_{e_1}(z_1-z_2\bo)\eo+P_2(\mu_0)G_{e_2}(z_1+z_2\bo)\et
\end{eqnarray*}
\begin{eqnarray*}
 &=& [{W}_{e_1}(z_1-z_2\bo)+ P_1(\lambda_0)F_{e_1}(z_1-z_2\bo)+P_1(\mu_0)G_{e_1}(z_1-z_2\bo)]\eo\\
 &+& [{W}_{e_2}(z_1+z_2\bo)+ P_2(\lambda_0)F_{e_2}(z_1+z_2\bo)+P_2(\mu_0)G_{e_2}(z_1+z_2\bo)]\et\\
 &=& [{w}_{e_1}(z_1-z_2\bo)]\eo+[{w}_{e_2}(z_1+z_2\bo)]\et.
\end{eqnarray*}
where $P_1({\lambda_0}),P_2({\lambda_0}),P_1({\mu_0}),P_2({\mu_0})\in\mathbb{R}$.$\Box$

\smallskip\smallskip
\noindent Now, using the last theorem and the decomposition:
$$\frac{w(z_1+z_2\bt) - \lambda _0 F(z_1+z_2\bt) - \mu _0 G(z_1+z_2\bt)}{(z_1+z_2\bt)-(z_{0,1}+z_{0,2}\bt)}$$
\begin{eqnarray*}
 &=& \frac{w_{e_1}(z_1-z_2\bo) - P_1(\lambda _0) F_{e_1}(z_1-z_2\bo) - P_1(\mu _0) G_{e_1}(z_1-z_2\bo)}{(z_1-z_2\bo)-(z_{0,1}-z_{0,2}\bo)}\eo\\
 &+& \frac{w_{e_2}(z_1+z_2\bo) - P_2(\lambda _0) F_{e_2}(z_1+z_2\bo) - P_2(\mu _0) G_{e_2}(z_1+z_2\bo)}{(z_1+z_2\bo)-(z_{0,1}+z_{0,2}\bo)}\et
\end{eqnarray*}
we obtain the following connections with the classical
theory of pseudoanalytic functions.

\begin{theorem}
Let $D_0=D_1\times_e D_2$ where $D_1$ and $D_2$ are open domains in $\mathbb{C}(\bo)$.
If $(F,G)$ is a $\bje$-generating pair in $D_{0}$ with
$w:D_{0}\subset\mathbb{T}\to\mathbb{T}$ a $(F,G)_{\bj}$-pseudoanalytic function on $D_{0}$ then
\begin{equation}
w(z_1+z_2\bt)=[{w}_{e_1}(z_1-z_2\bo)]\eo+[{w}_{e_2}(z_1+z_2\bo)]\et
\end{equation}
where ${w}_{e_k}$ is a $(F_{e_k},G_{e_k})$-pseudoanalytic function on $D_k$ for $k=1,2$.
Moreover,
\begin{equation}
\dot w(z_{1}+z_{2}\bt)=[\dot w_{e_1}(z_{1}-z_{2}\bo)]\eo+[\dot w_{e_2}(z_{1}+z_{2}\bo)]\et
\end{equation}
on $D_0$.
\label{central}
\end{theorem}

\begin{theorem}
If $w_{ek}:D_k\longrightarrow \mathbb{C}(\bold{i_1})$
is a $(F_{e_k},G_{e_k})$-pseudoanalytic function on the open domain $D_k$ for $k=1,2$
then the function $w:D_1\times_e D_2\longrightarrow \mathbb{T}$ defined as
$$w(z_1+z_2\bold{i_2})=w_{e1}(z_1-z_2\bold{i_1})\bold{e_1}+w_{e2}(z_1+z_2\bold{i_1})\bold{e_2} \mbox{ }\forall\mbox{ }z_1+z_2\bold{i_2}\in D_1\times_e D_2$$
is a $(F,G)_{\bj}$-pseudoanalytic function on $D_1\times_e D_2$ and
$$\dot w(z_{1}+z_{2}\bt)=\dot w_{e_1}(z_{1}-z_{2}\bo)\eo+\dot w_{e_2}(z_{1}+z_{2}\bo)\et$$
$\mbox{ }\forall\mbox{ }z_1+z_2\bold{i_2}\in D_1\times_e D_2.$
\label{central2}
\end{theorem}


\hspace{0.5cm} The last theorem gives another interpretation of Theorems \ref{central} and \ref{central2} in terms of Vekua equations.

\begin{theorem}
If $(F_{e_1},G_{e_1})$ and $(F_{e_2},G_{e_2})$ are complex (in \bo) generating pairs in
respectively $D_1$ and $D_2$. Then $w$
is a solution on $D_0=D_1\times_e D_2$ of the $\bj$-bicomplex Vekua equations with
the $\bje$-generating pair $(F,G)$ if and only if $w(z_1+z_2\bold{i_2})=w_{e1}(z_1-z_2\bold{i_1})\bold{e_1}+w_{e2}(z_1+z_2\bold{i_1})\bold{e_2}$
where $w_{e_k}$ is a solution on $D_k$ of the complex (in $\bo$) Vekua equation with the generating pair
$(F_{e_k},G_{e_k})$ for $k=1,2$.
\end{theorem}

\subsection{The Complexified Schr{\"o}dinger Equation}

Consider the equation
\begin{equation}
(\triangle_{\mathbb{C}}-\nu(z_1,z_2))f=0
\label{CL}
\end{equation}
in $\Omega\subset\mathbb{R}^4$, where $\triangle_{\mathbb{C}}={\partial^{2}_{z_1}}+{\partial^{2}_{z_2}}$,
$\nu$ and $f$ are complex (in $\bo$) valued functions. The equation (\ref{CL}) is simply the
complexification of the two-dimensio\-nal stationary Schr{\"o}dinger equation where
$\triangle_{\mathbb{C}}$ is the complex Laplacian (see \cite{28,23,27,26}).

\subsubsection{The Complex Laplacian}

First of all, we will write the complex Laplacian in a more explicit way to see that it contains
in the same time the classical Laplacian operator and the wave operator.

\begin{lemma}
Let $\omega=z_1+z_2\bt$, where $z_1,z_2\in \mathbb{C}(\bo)$ then
$$\partial_\omega \partial_{\bar{\omega}}=\frac{1}{4}({\partial^{2}_{z_1}}+{\partial^{2}_{z_2}})
=\frac{1}{4}\triangle_{\mathbb{C}}$$
$\forall f\in C^{2}(\Omega)$ where $\Omega$ is an open set in $\mathbb{R}^4$.
\end{lemma}
\emph{Proof.}
Let $\partial _{\omega } = \frac{1}{2}\left( {\partial _{z_1}  - \bt\partial _{z_2} } \right)$ and
$\partial _{\bar{\omega}} = \frac{1}{2}\left( {\partial _{z_1}  + \bt\partial _{z_2} } \right)$
then
\begin{eqnarray*}
4\partial_\omega \partial_{\bar{\omega}} &=& \partial _{z_1} \left( \partial _{z_1}  + \bt\partial _{z_2} \right) -
\bt\partial _{z_2}\left( \partial _{z_1}  + \bt\partial _{z_2} \right)\\
                                         &=& \partial^{2}_{z_1}+\bt\partial^{2}_{z_1 z_2}-\bt\partial^{2}_{z_2 z_1}+\partial^{2}_{z_2}\\
                                         &=& {\partial^{2}_{z_1}}+{\partial^{2}_{z_2}}.\Box
\end{eqnarray*}

\begin{proposition}
Let $\partial _{z_1 }  = \frac{1}{2}\left( {\partial _x  - \bo\partial _y }\right)$ and
$\partial _{z_2 }  = \frac{1}{2}\left( {\partial _p  - \bo\partial _q }\right)$ then
\begin{equation}
16\partial_\omega \partial_{\bar{\omega}} =4\triangle_{\mathbb{C}}=\left( {\partial^{2}_{x}}-{\partial^{2}_{y}}+{\partial^{2}_{p}}-{\partial^{2}_{q}}\right)
-2\bo\left( {\partial^{2}_{xy}}+{\partial^{2}_{pq}}  \right)
\end{equation}
$\forall f\in C^{2}(\Omega)$ where $\Omega$ is an open set in $\mathbb{R}^4$.

\end{proposition}
\emph{Proof.}
Consider,
\begin{eqnarray*}
4 {\partial^{2}_{z_1}} &=& \partial _{x}\left( {\partial _x  - \bo\partial _y } \right)-\bo\partial _{y}\left( {\partial _x  - \bo\partial _y } \right)\\
                       &=& \partial^{2}_{x}-\bo\partial^{2}_{xy}-\bo\partial^{2}_{yx}-\partial^{2}_{y}\\
                       &=& \partial^{2}_{x}-\partial^{2}_{y}-2\bo\partial^{2}_{xy}.
\end{eqnarray*}
Therefore,
\begin{eqnarray*}
4\left({\partial^{2}_{z_1}}+ {\partial^{2}_{z_2}}\right) &=& \left( {\partial^{2}_{x}}-{\partial^{2}_{y}}+{\partial^{2}_{p}}-{\partial^{2}_{q}}\right)
-2\bo\left( {\partial^{2}_{xy}}+{\partial^{2}_{pq}}  \right).\Box
\end{eqnarray*}

\smallskip
\begin{remark}
In the last proposition, if we let $y$ and $q$ be constant variables, then
\begin{enumerate}
\item $\Omega\subset\mathbb{C}(\bt)$;
\item $4\partial_\omega \partial_{\bar{\omega}}=\partial_z \partial_{\bar{z}}$
where $\partial _{\bar{z}} = \frac{1}{2}\left( {\partial _{x}  + \bt\partial _{p} }\right)
\mbox{ and } \partial _{z} = \frac{1}{2}\left( {\partial _{x}  - \bt\partial _{p} } \right)$;

\item $4\triangle_{\mathbb{C}}=4\partial_z \partial_{\bar{z}}={\partial^{2}_{x}}+{\partial^{2}_{p}}=\triangle, \mbox{ \textbf{the Laplacian operator}}.$

\end{enumerate}
Similarly, if $y$ and $p$ are constant variables, then
\begin{enumerate}
\item $\Omega\subset\mathbb{D}$;
\item $4\partial_\omega \partial_{\bar{\omega}}=\partial_z \partial_{\bar{z}}$
where $\partial _{\bar{z}} = \frac{1}{2}\left( {\partial _{x}  - \bj\partial _{q} }\right)
\mbox{ and } \partial _{z} = \frac{1}{2}\left( {\partial _{x}  + \bj\partial _{q} } \right)$;
\item $4\triangle_{\mathbb{C}}=4\partial_z \partial_{\bar{z}}={\partial^{2}_{x}}-{\partial^{2}_{q}}=\square, \mbox{ \textbf{the wave operator}}.$
\end{enumerate}
\label{both}
\end{remark}

\subsubsection{Factorization of the Complexified Schr{\"o}dinger Operator}

It is well known that if $f_0$ is a nonvanishing particular solution of the
one-dimensional stationary Schr{\"o}dinger equation
$$\left( -\frac{d^2}{dx^2}+\nu(x) \right)$$
then the Scr{\"o}dinger operator can be factorized as follows:
$$-\frac{d^2}{dx^2}+\nu(x)=\left(\frac{d}{dx}+\frac{f_0'}{f_0}\right)\left(\frac{d}{dx}-\frac{f_0'}{f_0}\right).$$

\smallskip\smallskip
The next result gives the analogue for the complexified Schr{\"o}dinger operator. By $C$ we denote the $\dagger_2$-bicomplex conjugation operator.

\begin{theorem}
Let $f_0:\Omega\subset\mathbb{R}^4\longrightarrow \mathbb{C}(\bo)$ be a nonvanishing particular solution of (\ref{CL}). Then for any
$\mathbb{C}(\bo)$-valued continuously twice differentiable function $\varphi$ the following equality hold:
\begin{equation}
(\triangle_{\mathbb{C}}-\nu)\varphi=4\left(\partial_{\bar{\omega}}+\frac{\partial_{\omega} f_0}{f_0}C\right)
\left(\partial_{\omega}-\frac{\partial_{\omega} f_0}{f_0}C\right)\varphi.
\end{equation}
\end{theorem}
\emph{Proof.}
Let $f_0:\Omega\subset\mathbb{R}^4\longrightarrow \mathbb{C}(\bo)$ be a nonvanishing particular solution of (\ref{CL}). Then
\footnotesize
\begin{eqnarray*}
\left(\partial_{\bar{\omega}}+\frac{\partial_{\omega} f_0}{f_0}C\right)
\left(\partial_{\omega}-\frac{\partial_{\omega} f_0}{f_0}C\right)\varphi &=& \left(\partial_{\bar{\omega}}+\frac{\partial_{\omega} f_0}{f_0}C\right)
\left(\partial_{\omega}\varphi-\frac{\partial_{\omega} f_0}{f_0}\varphi\right)\\
&=& \partial_{\bar{\omega}}\partial_{\omega}\varphi - \partial_{\bar{\omega}}\left(\frac{\partial_{\omega} f_0}{f_0}\varphi\right)
+\frac{\partial_{\omega} f_0}{f_0}\partial_{\bar{\omega}}\varphi-\frac{|\partial_{\omega} f_0|^{2}_{\bo}}{f^2_0}\varphi\\
&=& \frac{1}{4}\triangle_{\mathbb{C}}\varphi- \partial_{\bar{\omega}}\left(\frac{\partial_{\omega} f_0}{f_0}\right)
\varphi-\frac{|\partial_{\omega} f_0|^{2}_{\bo}}{f^2_0}\varphi\\
&=& \frac{1}{4}\triangle_{\mathbb{C}}\varphi- \left(\frac{\partial_{\bar{\omega}}\partial_{\omega} f_0\cdot f_0 -|\partial_{\omega} f_0|^{2}_{\bo}}{f^2_0}\right)\varphi
-\frac{|\partial_{\omega} f_0|^{2}_{\bo}}{f^2_0}\varphi\\
&=& \frac{1}{4}\triangle_{\mathbb{C}}\varphi-\frac{\frac{1}{4}\triangle_{\mathbb{C}}f_0}{f^2_0}\varphi.
\end{eqnarray*}
\normalsize
However, $(\triangle_{\mathbb{C}}-\nu)f_0=0 \Rightarrow \frac{\triangle_{\mathbb{C}}f_0}{f_0}=\nu$.
Hence,
\begin{equation}
(\triangle_{\mathbb{C}}-\nu)\varphi=4\left(\partial_{\bar{\omega}}+\frac{\partial_{\omega} f_0}{f_0}C\right)
\left(\partial_{\omega}-\frac{\partial_{\omega} f_0}{f_0}C\right)\varphi.\Box
\end{equation}

\begin{remark}
From the Remark \ref{both}, we see that the complexified Schr{\"o}dinger equation contains
the stationary two-dimensional Schr{\"o}dinger equation
$$(\triangle-\nu(x,p))f=0$$
and the Klein-Gordon equation
$$(\square-\nu(x,q))f=0.$$
Hence, our factorization of the complexified Schr{\"o}dinger equation is a generalization of
the factorization obtained in \cite{8} for the stationary two-dimensional Schr{\"o}dinger equation and
for the factorization obtained in \cite{24} for the Klein-Gordon equation.

\end{remark}

\subsubsection{Relationship Between Bicomplex Generalized Analytic Functions and Solutions of the Complexified
Schr{\"o}dinger Equation}

The next Lemma has been inspired from a similar result in the complex plane (see \cite{25}, p 140).

\begin{lemma}
Let $b:\Omega\subset\mathbb{R}^4\longrightarrow \mathbb{T}$ be a bicomplex function such that $b_{\omega}$
is $\mathbb{C}(\bo)$-valued, and let $W=u+\bt v:\Omega\subset\mathbb{R}^4\longrightarrow \mathbb{T}$ be a solution of the equation
\begin{equation}
W_{{\omega}^{\dagger_2}}=bW^{\dagger_2} \mbox{ on } \Omega.
\label{di1}
\end{equation}
Thus, $u:\Omega\subset\mathbb{R}^4\longrightarrow \mathbb{C}(\bo)$ is a solution of the equation
\begin{equation}
\partial_{{\omega}^{\dagger_2}}\partial_{\omega}u-(|b|^{2}_{\bo}+b_{\omega})u=0 \mbox{ on } \Omega
\label{di2}
\end{equation}
and $v:\Omega\subset\mathbb{R}^4\longrightarrow \mathbb{C}(\bo)$ is a solution of the equation
\begin{equation}
\partial_{{\omega}^{\dagger_2}}\partial_{\omega}v-(|b|^{2}_{\bo}-b_{\omega})v=0 \mbox{ on } \Omega.
\label{di3}
\end{equation}
\label{dd}
\end{lemma}
\emph{Proof.}
Using the $\dagger_2$ on both sides of Eq. (\ref{di1}), we have that
\begin{equation}
\partial_{{\omega}^{\dagger_2}}(u+\bt v)=b(u-\bt v) \Leftrightarrow \partial_{{\omega}}(u-\bt v)=b^{\dagger_2}(u+\bt v) \mbox{ on } \Omega.
\end{equation}
Therefore,
\begin{eqnarray*}
\partial_{\omega}\partial_{{\omega}^{\dagger_2}}(u+\bt v) &=& \partial_{\omega}b \cdot(u-\bt v)+b\partial_{\omega}(u-\bt v)\\
                                            &=& b_{\omega}(u-\bt v)+bb^{\dagger_2}(u+\bt v)\\
                                            &=& b_{\omega}(u-\bt v)+|b|^{2}_{\bo}(u+\bt v) \mbox{ on } \Omega.
\end{eqnarray*}
Now, by equality of the scalar and vectorial part, we obtain the Equations (\ref{di2}) and (\ref{di3}).$\Box$

\begin{theorem}
Let $W$ be a solution of the following bicomplex Vekua equation
\begin{equation}
\left( \partial_{{\omega}^{\dagger_2}} - \frac{\partial_{{\omega}^{\dagger_2}}f_0}{f_0}C \right) W=0
\label{gene}
\end{equation}
where $f_0$ is a nonvanishing solution of the complexified Schr{\"o}dinger equation (\ref{CL}).
Then $u=\mbox{Sc}(W)$ is a solution of (\ref{CL}) and $v=\mbox{Vec}(W)$ is a solution of
the equation
\begin{equation}
\left(\triangle_{\mathbb{C}}+\nu(z_1,z_2)-2{\left( \frac{|\bigtriangledown_{\mathbb{C}}f_0|_{\bo}}{f_0} \right)}^{2}\right)v=0
\end{equation}
where $\bigtriangledown_{\mathbb{C}}={\partial_{z_1}}+\bt {\partial_{z_2}}$.
\end{theorem}
\emph{Proof.}
Consider the function $b=\frac{\partial_{{\omega}^{\dagger_2}}f_0}{f_0}$. Then
\begin{eqnarray*}
b_{\omega}=\partial_{\omega}\left(\frac{\partial_{{\omega}^{\dagger_2}}f_0}{f_0}\right)
&=& \frac{(\partial_{\omega}\partial_{{\omega}^{\dagger_2}}f_0)\cdot f_0-(\partial_{{\omega}^{\dagger_2}}f_0)(\partial_{\omega}f_0)}{f^2_2}\\
&=& \frac{\triangle_{\mathbb{C}}f_0}{4f_0}-\frac{|\partial_{{\omega}^{\dagger_2}}f_0|^{2}_{\bo}}{f^2_2}\\
&=& \frac{\nu(z_1,z_2)}{4}-\frac{|\partial_{{\omega}^{\dagger_2}}f_0|^{2}_{\bo}}{f^2_2}.
\end{eqnarray*}
Therefore, $b_{\omega}$ is a $\mathbb{C}(\bo)$-valued function.
Now, from Lemma \ref{dd},
\begin{equation*}
\frac{\triangle_{\mathbb{C}}u}{4}=\left( \frac{\nu(z_1,z_2)}{4}- \frac{|\partial_{{\omega}^{\dagger_2}}f_0|^{2}_{\bo}}{f^2_2}+ \frac{|\partial_{{\omega}^{\dagger_2}}f_0|^{2}_{\bo}}{f^2_2} \right)u
\end{equation*}
i.e $$(\triangle_{\mathbb{C}}-\nu)u=0,$$
and
\begin{eqnarray*}
\frac{\triangle_{\mathbb{C}}v}{4} &=& \left(-\frac{\nu(z_1,z_2)}{4}+\frac{|\partial_{{\omega}^{\dagger_2}}f_0|^{2}_{\bo}}{f^2_2}+ \frac{|\partial_{{\omega}^{\dagger_2}}f_0|^{2}_{\bo}}{f^2_2} \right)v\\
                                  &=& \left(-\frac{\nu(z_1,z_2)}{4}+\frac{2|\partial_{{\omega}^{\dagger_2}}f_0|^{2}_{\bo}}{f^2_2}\right)v\\
                                  &=& \left(-\frac{\nu(z_1,z_2)}{4}+\frac{(\partial_{z_1}f_0)^2+(\partial_{z_2}f_0)^2}{2f^2_2}\right)v
\end{eqnarray*}
i.e. $$\left( \triangle_{\mathbb{C}}-\eta \right)v=0$$
where $\eta(z_1,z_2)=-\nu(z_1,z_2)+\frac{2|\bigtriangledown_{\mathbb{C}}f_0|^{2}_{\bo}}{f^2_2}.\Box$

\begin{remark}
From Theorem \ref{iifR1}, if $W$ possesses a $(f_0,\frac{\bt}{f_0})_{\bo}$-derivative on an open set $\Omega\subset\mathbb{T}$ then
$W$ is a solution of the bicomplex Vekua equation (\ref{gene}):
$$\left( \partial_{{\omega}^{\dagger_2}} - \frac{\partial_{{\omega}^{\dagger_2}}f_0}{f_0}C \right) W=0\mbox{ on }\Omega.$$
In that case, $a_{(F,G)}^{(2)}=0$ and $b_{(F,G)}^{(2)}=\frac{\partial_{{\omega}^{\dagger_2}}f_0}{f_0}$ where
$$F=f_0\quad \mbox{and}\quad G=\frac{\bt}{f_0}$$
is a ${\bo}$-generating pair for (\ref{gene}).
\end{remark}

From the last remark, we can conclude that the bicomplex pseudoanalytic function theory open the way to find explicit solutions of the
complexified Schr{\"o}dinger equation.


\newpage
\begin{thebibliography}{99}

\bibitem{1} S. Agmon and L. Bers, The expansion theorem for pseudo-analytic functions, {\em Proc. Amer. Math. Soc.} \textbf{3},
757--764 (1952).

\bibitem{25} K. W. Bauer, {\em Differential Operators for Partial Differential Equations and Function Theoretic Applications},
Springer, Berlin (1980).

\bibitem{2} H. Bahlouli, A. D. Alhaidari, A. Al-Zahrani and N. Economou, Study of electromagnetic wave propagation
in active medium and the equivalent Schrodinger equation with energy-dependant complex potential,
{\em Phys. Rev. B} \textbf{72} 94304 (2005).

\bibitem{3} L. Bers, An outline of the theory of pseudoanalytic function, {\em Bull. Am. Math. Soc.} \textbf{62},
291--331 (1956).

\bibitem{4} L. Bers, {\em Theory of Pseudo-Analytic Functions}, New York University, New York (1952).

\bibitem{5} A. Casta\~neda and V.V. Kravchenko, New applications of pseudoanalytic function theory to
the Dirac equation, {\em J. Phys. A.: Math. Gen.} \textbf{38}, No. 42,
9207--9219 (2005).

\bibitem{6} N. Fleury, M. Rausch de Traubenberg and R.M. Yamaleev, Commutative extended complex numbers and
connected trigonometry,  {\em J. Math. Ann. and Appl.} \textbf{180},
431--457 (1993).

\bibitem{28} K. Fujita and M. Morimoto, Conical Fourier transform of Hardy space of harmonic functions on the Lie ball,
{\em Annales Universitatis Mariae Curie - Sklodowska} \textbf{LIII}, No. 4, 41--55 (1999).

\bibitem{23} K. Fujita and M. Morimoto, Integral representation for eigenfunctions of the Laplacian,
{\em J. Math. Soc. Japan} \textbf{51}, No. 3, 699--713 (1999).

\bibitem{27} K. Fujita and M. Morimoto, Reproducing Kernels Related to the Complex Sphere,
{\em Tokyo J. Math.} \textbf{23}, No. 1, 161--185 (2000).

\bibitem{7} P.E. Hodgson, {\em The Nuclear Optical Model: Introductory Overview}, Bruyeres-le-Chatel, France (1996).

\bibitem{8} V.V. Kravchenko, On a relation of pseudoanalytic function theory to the two-dimensional stationary
Schr{\"o}dinger equation and Taylor series in formal powers for its solutions, {\em J. Phys. A.: Math. Gen.} \textbf{38},
No. 18, 3947--3964 (2005).

\bibitem{20} V.V. Kravchenko, On a factorization of second-order elliptic operators and applications,
 {\em J. Phys. A} \textbf{39}, No. 40, 12407--12425 (2006).

\bibitem{9} V.V. Kravchenko, On the relationship between $p$-analytic functions and the Schr{\"o}dinger equation,
{\em Zeitschrift f{\"u}r Analysis und ihre Anwendungen} \textbf{24}, No. 3, 487--496 (2005).

\bibitem{24} V.V. Kravchenko, D. Rochon and S. Tremblay, On the Klein-Gordon Equation and the Hyperbolic Pseudoanalytic Function Theory,
(to appear).

\bibitem{26} M. Morimoto and K. Fujita, Analytic functionals and entire functionals on the complex light cone,
{\em Hiroshima Math. J.} \textbf{25}, 493--512 (1995).

\bibitem{10} N.F. Mott and H.S. Massey, {\em The theory of atomic collision}, Claredon Press, Oxford (1965).

\bibitem{11} G.B. Price, {\em An introduction to multicomplex spaces and functions},
Marcel Dekker Inc., New York (1991).

\bibitem{12} D. Rochon (2004), A bicomplex Riemann zeta function, {\em Tokyo J. Math.} \textbf{27}, 357--369 (2004).

\bibitem{13} D. Rochon, A generalized Mandelbrot set for bicomplex numbers, {\em Fractal} \textbf{8}, 355--368 (2000).

\bibitem{22} D. Rochon, Sur une g\'en\'eralisation des nombres complexes: les t\'etranombres,
(M. Sc. Universit\'e de Montr\'eal, 1997).

\bibitem{14} D. Rochon and M. Shapiro, On algebraic properties of bicomplex and hyperbolic numbers,
{\em Anal. Univ. Oradea}, fasc. math., {\bf 11}, 71--110 (2004).

\bibitem{15} D. Rochon and S. Tremblay, Bicomplex Quantum Mechanics: I. The Generalized Schr{\"o}dinger Equation,
{\em Adv. App. Cliff. Alg.} {\bf 12}, No. 2, 231--248 (2004).

\bibitem{19} D. Rochon and S. Tremblay, Bicomplex Quantum Mechanics: II. The Hilbert Space,
{\em Adv. App. Cliff. Alg.} {\bf 16}, No. 2, 135--157 (2006).

\bibitem{21} B. V. Shabat, {\em Introduction to Complex Analysis part II: Functions of Several Variables},
(American Mathematical Society, 1992).

\bibitem{16} G. Sobczyk, The hyperbolic number plane, {\em Coll. Maths. Jour.} \textbf{26}, No. 4, 268--280 (1995).

\bibitem{17} I. N. Vekua, {\em Generalized Analytic Functions}, Pergamon Press Ltd (English translation), Oxford (1962).

\bibitem{18} Guo-chun Wen, {Linear and Quasilinear Complex Equations of Hyperbolic and Mixed Type}, Taylor \& Francis, London (2002).


\end{thebibliography}
\end{document}